\documentstyle[12pt]{amsart}
\def\bP{{\Bbb P}}

\def\C{{\Bbb C}}
\def\Z{{\Bbb Z}}

\def\cD{{\cal D}}
\def\cV{{\cal V}}
\def\cL{{\cal L}}
\def\O{{\cal O}}
\def\si{\sigma }
\def\bC{{\bar C}}
\def\a{\alpha }
\def\b{\beta }
\def\om{\omega }

\def\p{\phi }
\def\r{\rho }
\def\th{\theta }
\def\ga{\gamma }

\begin{document}
\baselineskip=16.5pt

\title[Jet bundles over the projective line]
{A remark on the jet bundles over the projective line}
\author{Indranil Biswas}

\address{School of Mathematics, Tata Institute of Fundamental
Research, Homi Bhabha Road, Bombay 400005, INDIA}
\email{indranil@@math.tifr.res.in}
\date{}

\maketitle

\section{Introduction}

Let $X$ be a Riemann surface equipped with a projective
structure (i.e., a covering by coordinate charts such that
the transition functions are of the form $z \longmapsto
(az +b)/(cz+d)$). Let $\cL$ be a line bundle on $X$ such that
${\cL}^2 = T_X$. Let $J^m({\cL}^{\otimes n}) \longrightarrow X$
denote the jet bundle of
order $m$ for the line bundle ${\cL}^{\otimes n}$. For $i \geq j$,
there is a natural restriction homomorphism from
$J^i({\cL}^{n})$ onto $J^j({\cL}^{n})$. We prove that
for any $m\geq n$, the surjective homomorphism
$$
J^m({\cL}^{n}) \, \longrightarrow \, J^n({\cL}^{n})
$$
admits a canonical splitting [Theorem 4.1]. As a consequence,
for each $n\geq 0$ we
construct a differential operator of order
$n$ from ${\cL}^{n -1}$ to ${\cL}^{-n-1}$ whose symbol is the
constant function $1$. Theorem 4.1 follows from the results on the
jet bundles over the projective line established in Section 2.

In \cite{CMZ} certain
differential operators on a Riemann surface
equipped with a projective structure are explicitly constructed
(see (3.1)). As an application
of the set-up we use to prove Theorem 4.1, in Section 3 we derive the 
differential operators
constructed in \cite{CMZ}. The present work was inspired by
\cite{CMZ}; in fact, it grew out of attempts to
reconstruct the differential operators there without using
the coordinates.

\section{Constructions on the projective line}

Let $V$ be a two dimensional vector space over $\C$. Let $\bP(V)$
denote the projective space given by the space of all one
dimensional quotients of $V$. Define the line bundle
$$
L \, := \, {\O}_{\bP(V)}(1)
$$
on $\bP(V)$, whose fiber over the quotient line $[q]$ is
the line $[q]$ itself.

We will recall the definition of the jet bundles for a line bundle.
For a line bundle $\xi$ on a Riemann surface $X$, the $n$-th
jet bundle, denoted by $J^n(\xi)$, is the rank $n+1$ vector
bundle on $X$ whose fiber over $x \in X$ is
$$
{\xi}_x {\otimes}_{\C} ({\O}_{X,x}/{\bf m}^{n+1}_x)
$$
where ${\xi}_x$ is the fiber of $\xi$ over $x$, and ${\O}_{X,x}$ is the
ring of functions defined around $x$ with ${\bf m}_x$ being the
maximal ideal consisting of functions vanishing at $x$.
The inclusion of ${\bf m}^{n+1}_x$ in
${\bf m}^n_x$ induces the following short exact
sequence of vector bundle on $X$:
$$
0 \, \longrightarrow  \, K^{\otimes n}_X\otimes \xi \, 
\longrightarrow \, J^n(\xi) \, \longrightarrow \, J^{n-1}(\xi) \, 
\longrightarrow \, 0 \leqno{(2.1)}
$$
where $K_X$ is the canonical bundle of $X$.

Let $\cV$ denote the rank two trivial vector bundle on $\bP(V)$ with
$V$ as the fiber. For any $n \geq 0$, $S^n(\cV)$ will denote the $n$-th
symmetric power of $\cV$, with $S^0(\cV)$ being the trivial line
bundle.

\medskip
\noindent {\bf Lemma 2.2.}\, {\it For any integer $n \geq 0$, the 
vector bundle $J^n(L^n)$ on
$\bP(V)$ is canonically isomorphic to the symmetric power $S^n(\cV)$.
For any $m \geq n$, the surjection
$$
J^m(L^n) \,  \longrightarrow \, J^n(L^n) \, \longrightarrow \, 0
$$
given by (2.1), admits a canonical splitting (i.e., a homomorphism
from $J^n(L^n)$ to $J^m(L^n)$ such that the composition is identity
on $J^n(L^n)$).}
\medskip

\noindent {\bf Proof.}\, Take any integer $n \geq 0$. Since $S^n(V) =
H^0(\bP(V), L^n)$, for any $x \in \bP(V)$ there is a natural homomorphism
of $S^n(V)$ into the
fiber $J^n(L^n)_x$ given by the restriction of sections to the $n$-th
order infinitesimal neighborhood of $x$. Since for any integer $j$ with
$j\leq n$,
$$
\dim H^0(\bP(V), L^n \otimes {\O}_{\bP(V)}(-jx)) \, - \,
\dim H^0(\bP(V), L^n\otimes {\O}_{\bP(V)}(-(j+1)x)) \, = \, 1
$$
the above obtained homomorphism must be an isomorphism.
This proves the first part of the lemma.

Take any integer $m$ such that $m\geq n$. We may restrict a section
of $L^n$ to the $m$-th order infinitesimal neighborhood
of $x$ to get a homomorphism from the vector space $S^n(V)$ (=
$H^0(\bP(V), L^n)$) to the fiber $J^m(L^n)_x$. Now
using the previous identification of $S^n(V)$ with $J^n(L^n)_x$ we
get the required splitting. $\hfill{\Box}$
\medskip

Setting $m = n+1$ in Lemma 2.2 we obtain the following:

\medskip
\noindent {\bf Corollary 2.3.}\, {\it For any integer
$n \geq 0$, the exact sequence
$$
0\, \longrightarrow \, K^{n+1}_{\bP(V)}\otimes L^n \, \longrightarrow
 \, J^{n+1}(L^n) \, \longrightarrow \, J^n(L^n)\, \longrightarrow \,0
$$
admits a canonical splitting.}
\medskip

Choose and fix a trivialization of ${\stackrel{2}{\wedge}V}$; this is
equivalent to fixing a nonzero vector $\th$ in ${\stackrel{2}{\wedge}}V$.
The canonical bundle $K_{\bP(V)} = L^{-2}\otimes {\rm det}\, \cV$. Using
the trivialization of ${\stackrel{2}{\wedge}}V$ we have,
$K_{\bP(V)} = L^{-2}$.

The sheaf of differential operators of order $k$ from the sections
of a line bundle $\xi$ to the sections of a line bundle $\eta$
is precisely the sheaf ${\rm Hom}(J^k(\xi), \eta)$. Consider
the projection of $J^{n+1}(L^n)$ onto
$K^{n+1}_{\bP(V)}\otimes L^n =  L^{-n-2}$ defining the splitting in 
Corollary 2.3. This gives a global differential operator of order $n+1$,
$$
{\cD}({n+1}) \, \in \, H^0(\bP(V), {\rm Diff}^{n+1}(L^n,
L^{-n-2})) \leqno{(2.4)}
$$

The symbol of a differential operator in ${\rm Diff}^{n+1}(L^n, L^{-n-2})$
is a section of of the line bundle $T^{n+1}_{\bP(V)}\otimes L^{-2n-2} = 
\O_{\bP(V)}$. 
Since the differential operator ${\cD}(n+1)$ in (2.4) gives a splitting of
the jet sequence -- it's symbol, which is a constant function, must be
the constant function $1$.

Let $SL(V)$ denote the subgroup of $GL(V) = {\rm Aut}(V)$ that acts
trivially
on ${\stackrel{2}{\wedge}}V$. The group $SL(V)$ has a natural action
on $\bP(V)$, and ${\rm Aut}(\bP(V)) = SL(V)/{\Z}_2$. There is a
natural induced
action of $SL(V)$ on any sheaf $J^m(L^n)$ that lifts the action
on $\bP(V)$. The isomorphism between $S^n(\cV)$
and $J^n(L^n)$, and the splitting in Lemma 2.2, are both equivariant for
this action. Indeed, this follows from the
canonical nature of the construction in Lemma 2.2. So, in particular,
the differential operator ${\cD}(n)$ in (2.4) is an invariant for the
action of $SL(V)$ on the space of all global sections of ${\rm 
Diff}^n(L^{n-1}, L^{-n-1})$.

Note that ${\cD}(n)$ is not an invariant for the action $GL(V)$ since
the trivialization of ${\stackrel{2}{\wedge}}V$ was used in its
construction. The identification between $K_{\bP(V)}$ and $L^{-2}$
is not equivariant for the action of the center of $GL(V)$.

Setting $n=2$ in Lemma 2.2 we get that $J^2(L^2) = S^2(\cV)$. This
implies that the homomorphism
$$
\r \, : \, H^0(\bP(V), J^2(L^2)) \, \longrightarrow \, H^0(\bP(V), L^2)
$$
induced by the obvious projection, namely $J^2(L^2) \longrightarrow L^2$,
is actually an isomorphism. Moreover, $\r$ is the identity map of
$S^2(V)$. Thus, after identifying the
tangent bundle $T_{\bP(V)}$ with $L^2$ using the trivialization
of ${\stackrel{2}{\wedge}}V$, the Lemma 2.2
implies that
$$
H^0(\bP(V) , T_{\bP(V)}) \, = \, H^0(\bP(V), J^2(L^2))
\, = \, S^2(V) \leqno{(2.5)}
$$
The Lie-bracket operation equips the vector space
$H^0(\bP(V), T_{\bP(V)})$ with the structure of a Lie algebra. The
action of $SL(V)$ on $\bP(V)$ gives a Lie algebra homomorphism from
its Lie algebra, $sl(V)$, into $H^0(\bP(V), T_{\bP(V)})$. This
homomorphism is actually an isomorphism. The Lie
algebra structure on $S^2(V)$ induced by the equality (2.5) can be
seen directly as follows: using contraction, $S^2(V)$ maps
$V^*$ into $V$; on the other hand, $\th$ identifies $V^*$ with $V$ --
combining these, the resulting homomorphism from $S^2(V)$ into $sl(V)$
is an isomorphism.

Let $C \in S^2(H^0(\bP(V), T_{\bP(V)}))$ be the Casimir of the Lie
algebra $H^0(\bP(V), T_{\bP(V)})$. The section $C$ is evidently
an invariant for the obvious action of $SL(V)$ on the vector
space $S^2(H^0(\bP(V), T_{\bP(V)}))$. For a section $s$ of 
$T_{\bP(V)}$, let $L_{s}$ denote
the Lie derivative with respect to $s$. The (second order) Lie
derivative with respect to $s{\otimes} s$ is defined to be
$L_{s}\circ L_{s}$. Thus $C$ acts as a differential operator,
denoted by $L_C$, on all vector bundles associated to $\bP(V)$. This
differential operator is actually of order zero (i.e., a constant
scalar multiplication).

Using (2.5) and Lemma 2.2 we get that
$$
S^2(H^0(\bP(V), T_{\bP(V)})) \, = \, H^0(\bP(V) ,S^2(J^2(T_{\bP(V)})))
$$
Let $\bC \in H^0(\bP(V), S^2(J^2(T_{\bP(V)})))$ be the element
corresponding to the Casimir $C$; $\bC$ is actually the Casimir for
the Lie algebra $S^2(V)$ (which is the fiber of $J^2(T_{\bP(V)})$).

Let $p : J^2(T_{\bP(V)}) \longrightarrow T_{\bP(V)}$ be the obvious
projection. If (locally)
$$
\bC \, = \, \sum_{i} A_i\otimes A_i \leqno{(2.6)}
$$
where $A_i$ are local sections of $J^2(T_{\bP(V)})$,
consider the operator
$$
L_{\bC} \, = \, \sum L_{p(A_i)}\circ L_{p(A_i)}
$$
with $L_{\p(A_i)}$ being the Lie derivative with respect to the
vector field $p(A_i)$. It is easy to 
check that the operator $L_{\bC}$ does not depend upon the choice of the
decomposition of $\bC$, and that $L_{C} = L_{\bC}$.

\section{Jets of the trivial line bundle on the projective line}

Let ${\rm Diff}^n(\O , \O) = J^n(\O)^*$ be the sheaf of differential
operators on the trivial line bundle over $\bP(V)$. The symbol map, which
is the dual of the injection in (2.1), gives a surjective homomorphism
$$
\si \, : \,  {\rm Diff}^n(\O, \O) \, \longrightarrow \, T^n_{\bP(V)} 
$$
Let $\ga$ denote the obvious projection of $J^{n-1}(T^n_{\bP(V)})$ onto
$T^n_{\bP(V)}$.

For any $n \geq 1$, let $J^n_0(\O) \subset J^n(\O)$ be the kernel of
the obvious homomorphism from $J^n(\O)$ onto $J^0(\O) = \O$. This subsheaf
has a canonical splitting given by the constant functions. Define the 
subsheaf, ${\rm Diff}^n_0(\O, \O) := J^n_0(\O)^*$,
of ${\rm Diff}^n(\O,\O)$.

For a function $f$ on $\C$ and any integer $n \geq 1$, in Proposition 1 
(page 4) of \cite{CMZ} (where it is called ${\cL}_{-n}(f)$) the following 
differential operator of order $n$ is constructed:
$$
D_{n}(f) \, := \, \sum_ {i=0}^{n-1} {{(2n-i)!}\over {i!(n-i)!(n-i-1)!}}
f^{(i)}{\partial}^{n-i}  \leqno{(3.1)}
$$
with $f^{(i)} = {\partial}^{i}f$ being the $i$-th derivative of $f$. 
The operator
${D}_{n}$ has the property that for any M\"obius transformation,
$M(z) = (az+b)/(cz+d)$, of $\C{\Bbb P}^1$, the following equality holds:
$$
{D}_{n}(f)\circ M \, = \, {D}_{n}(M_n.(f \circ M)) \leqno{(3.2)}
$$
where $M_n(z) =(cz+d)^{2n}$.

Thus ${D}_{n}$ is a $SL(V)$ 
equivariant ${\O}_{\bP(V)}$ linear isomorphism (in other words, a
canonical isomorphism)
$$
\p \, : \,  J^{n-1}(T^n_{\bP(V)}) \, \longrightarrow \, 
{\rm Diff}^n_0(\O,\O) \leqno{(3.3)}
$$
The operator ${D}_{n}$ has the further property that $\si \circ \p = 
\ga$. It is shown in \cite{CMZ} that this splitting condition together
with the automorphic property (3.2) actually determine the operator 
${D}_{n}$. In this section we 
want to deduce the above result of \cite{CMZ} in the set-up of Section 2.

Take a point $x \in \bP(V)$. The long exact sequence of cohomology
for the exact sequence of sheaves on $\bP(V)$
$$
0 \, \longrightarrow \, {\O}(-(n+1).x) \, \longrightarrow \, \O
\, \longrightarrow \, J^n(\O)_x \, \longrightarrow \, 0
$$
gives the equality
$$
J^n_0(\O)_x \, = \,  H^1(\bP(V) , {\O}(-(n+1).x))
$$
where $J^n_0(\O)_x$ is the fiber of $J^n_0(\O)$ over $x$.

Choose and fix an isomorphism between the two line bundles ${\O}(x)$
and $L$. Since the fiber ${\O}(x)_x = T_{\bP(V), x} = L^2_x$, fixing
such an isomorphism is equivalent
to fixing a nonzero vector $\om$ in $L_x$.

Using Serre duality for ${\O}(-(n+1).x)$, and then
identifying $K_{\bP(V)}$ with ${\O}(-2x)$ using $\om$, we have
$$
{\rm Diff}^n_0(\O, \O)_x \, = \, H^0(\bP(V), \O((n-1).x))
\, = \, H^0(\bP(V), L^{n-1}) \, = \, S^{n-1}(V) \leqno{(3.4)}
$$

Consider the restriction of sections of $T^n_{\bP(V)}$ to the
$(n-1)$-th order infinitesimal neighborhood of $x$, namely
$$
\b \, : \, S^{2n}(V) \, = \, H^0(\bP(V), T^n_{\bP(V)})
\, \longrightarrow \, J^{n-1}(T^n_{\bP(V)})_x \leqno{(3.5)}
$$
which is clearly a surjective homomorphism. Indeed, in the proof
of Lemma 2.2 we saw that $S^{2n}(V)$ surjects onto $J^{2n}(L^{2n})_x
= J^{2n}(T^n_{\bP(V)})_x$. We want to identify the kernel of the
homomorphism $\b$.

The symplectic form on $V$ given by the trivialization of
${\stackrel{2}{\wedge}}V$ identifies $V$ with $V^*$.
Let $v$ be the vector in the kernel of the quotient homomorphism
$V \longrightarrow L_x$ which corresponds to $\om$ using the symplectic
form on $V$. (This vector $v \in V$ corresponds to the section
of the sheaf ${\O}(x)$ given by the constant function $1$.)

Consider the homomorphism, $m_v : S^n(V) \longrightarrow S^{2n}(V)$,
defined by multiplication with $v^{\otimes n}$. The inclusion
$m_v$ corresponds to the natural inclusion of $H^0(\bP(V), {\O}(n.x))$
into $H^0(\bP(V), \O(2n.x))$. The image of $m_v$ is precisely the kernel
of $\b$ in (3.5).

Consider the homomorphism $i_{\om} : S^{2n}(V) \longrightarrow S^{n-1}(V)$
given by the contraction with ${\om}^{\otimes (n+1)}$. (The vector $\om$
is considered as an element of $V^*$.) This homomorphism
vanishes on the image $m_v(S^n(V))$. Indeed, this follows from the
fact that $\om (v) = 0$.

Thus using the equality (3.4) and $i_{\om}$ we have the homomorphism
$$
{\p}_x \, : \, J^{n-1}(T^n_{\bP(V)})_x \, \longrightarrow \,
{\rm Diff}^n_0(\O, \O)_x
$$
It is easy to check that the homomorphism ${\p}_x$ does not depend
upon the choice of the nonzero vector $\om \in L_x$. The resulting 
homomorphism
$\p$ from $J^{n-1}(T^n_{\bP(V)})$ to ${\rm Diff}^n_0(\O,\O)$ satisfies
the condition that $\si \circ \p = \ga$. The canonical nature of the
construction of $\p$ ensures that it is equivariant for the action
of $SL(V)$. Since ${\p}_x$ is an isomorphism, $\p$ is an isomorphism.

Let $U \subset SL(V)$ be the unipotent subgroup which fixes the 
vector $v$. Let
${\frak n}$ be the nilpotent part of the Lie algebra of $U$. Let $N$
denote the unique element in $\frak n$ which maps a preimage of $\om$
(in $V$) to $v$. For any $0\leq i \leq 2n$, the image of $S^i(V)$ in
$S^{2n}(V)$, for the homomorphism given by the multiplication with 
$v^{2n-i}$, is denoted by
$S^{2n}_i(V)$. In this notation, $N$ maps $S^{2n}_{i+1}(V)$ onto 
$S^{2n}_i(V)$;
the resulting homomorphism from $S^{i+1}(V)$ onto $S^i(V)$ is
the contraction by $\om$. From this
it is easy to deduce that any homomorphism from $S^{2n}(V)/m_v(S^n(V))$
to $S^{n-1}(V)$, which is equivariant
for the actions of $N$, must be a scalar multiple of $i_{\om}$.
Now the condition, $\si \circ \p = \ga$, uniquely determines
the homomorphism $\p$.

\section{Jets on a Riemann surface with a projective structure}

Let $X$ be a Riemann surface, not necessarily compact.
A {\it projective structure}
on $X$ is a maximal atlas of holomorphic coordinate charts, $\{U_{\a},
f_{\a}\}_{\a\in I}$, covering $X$, such that any $f_{\a}$ maps $U_{\a}$
biholomorphically onto some analytic open set in $\bP(V)$ and the
transition function $f_{\a}\circ f^{-1}_{\b}$, for any $\a , \b \in I$,
is a restriction of an automorphism of $\bP(V)$, \cite{G}, \cite{D}, 
\cite{T}. We note that any Riemann surface admits a
projective structure, since, from the uniformization theorem, the
universal cover has a natural projective structure. It is know that for 
any projective structure, it is possible
to choose a sub-cover such that the transition functions have a compatible
lift to $SL(V)$ (from ${\rm Aut}(\bP(V))$). Actually, more than
one inequivalent lifts are possible. For a compact Riemann surface,
the set of equivalence classes of lifts correspond to the set of square 
roots of the canonical bundle (called
{\it theta characteristics}) \cite{G},
\cite{T}. Henceforth, by a projective structure we will always mean
a lift of the structure group to $SL(V)$.

Let $X$ be a Riemann surface equipped with
a projective structure in the above sense.

Since the natural action of $SL(V)$ on $\bP(V)$ lifts to the bundle $L$,
the projective structure gives a line bundle on $X$ associated to $L$.
Let $\cL$ denote this line bundle on $X$. Since the isomorphism
between $L^2$ and $T_{\bP(V)}$ is $SL(V)$ equivariant, we have
${\cL}^{\otimes 2} = T_X$.

Since $J^n(L^n)$ on $\bP(V)$ is a trivial bundle (Lemma 2.2), it has
a natural flat connection, which is equivariant under the action of
$SL(V)$. We now
have the following consequence of Lemma 2.2, Corollary 2.3 and (2.4):

\medskip
\noindent {\bf Theorem 4.1.}\, {\it For any $n \geq 0$, the jet bundle 
$J^n({\cL}^n)$ on
$X$ has a natural flat connection, and $S^n(J^1(\cL)) = J^n({\cL}^n)$,
with the identification being compatible with the flat connections.
For any $m\geq n$, the natural surjection
$$
J^m({\cL}^n) \, \longrightarrow \, J^n({\cL}^n) \, \longrightarrow \, 0
$$
has a canonical splitting. Setting $m=n+1$, a global differential
operator of order $n+1$
$$
{\cD}_X (n+1) \, \in \, H^0(X, {\rm Diff}^n_X({\cL}^{n}, {\cL}^{-n-2}))
$$
is obtained. The symbol of ${\cD}_X(n)$ is the constant function $1$.
The fibers of $J^2(T_X)$ have the structure of a Lie algebra compatible
with the flat connection on $J^2(T_X)$. The Lie derivative action
of the Casimir
$$
C_X \, \in \, H^0(X, S^2(J^2(T_X)))
$$
on any tensor power of $\cL$ is a multiplication by a constant scalar.}
\medskip

Similarly, since the isomorphism $\p$ in (3.3) is $SL(V)$ equivariant,
we have an isomorphism of vector bundles on $X$
$$
{\p}_X \, : \, J^{n-1}(T^n_X) \, \longrightarrow \, {\rm Diff}^n_0(\O, \O)
$$
(${\rm Diff}^n_0 (\O,\O) \subset {\rm Diff}^n_X(\O,\O)$ is the canonical
complement of ${\rm Diff}^0_X(\O,\O)$)  
such that the composition of the symbol map on ${\rm Diff}^n_X(\O, \O)$
with the isomorphism ${\p}_X$ is the natural projection
of $J^{n-1}(T^n_X)$ onto $T^n_X$.


\end{document}